\documentclass[12pt]{article}

\usepackage{amssymb}

\newtheorem{thm}{Theorem}
\newtheorem{rmk}{Remark}
\newtheorem{lem}{Lemma}
\newtheorem{defin}{Definition}

\begin{document}

\noindent \textbf{\large CONVERGENCE RATE OF WAVELET EXPANSIONS OF\\ GAUSSIAN RANDOM PROCESSES\footnotetext{Resubmitted: \today}}
\vskip 5mm

{\it Short title:} \textbf{WAVELET EXPANSIONS OF GAUSSIAN PROCESSES}

\vskip 1cm
\noindent \textbf{\large Yuriy Kozachenko$^a$, Andriy Olenko$^b$\footnote{Address correspondence to Andriy Olenko, Department of Mathematics and Statistics, La Trobe University, Victoria 3086, Australia; E-mail:  a.olenko@latrobe.edu.au} and Olga Polosmak$^c$}

\vskip 5mm
\noindent {$^a$ Department of Probability Theory, Statistics and Actuarial Mathematics, Taras Shevchenko Kyiv National University, Kyiv, Ukraine}

\noindent {$^b$ Department of Mathematics and Statistics, La Trobe University, Melbourne, Australia}

\noindent {$^c$ Department of Economic Cybernetics, Taras Shevchenko Kyiv National University, Kyiv, Ukraine}

\vskip 3cm

\vskip 3mm

\vskip 6mm

\noindent {\bf ABSTRACT}

\noindent \textit{The paper characterizes uniform convergence rate for general classes of wavelet expansions of
stationary Gaussian random processes. The convergence in probability is considered.}
\vskip 4mm

\newpage


\section{Introduction}
In various statistical, data compression, signal processing applications and simulation, it could be used to convert the problem of analyzing a continuous-time random process to that of analyzing a random sequence, which is much simpler. Multi\-resolution analysis provides an efficient framework for the decomposition of random processes. This approach is widely used in statistics to estimate a curve given observations of the curve plus some noise.

Various extensions of the standard statistical methodology were proposed recently. These include curve estimation in the presence of correlated noise. For these purposes the wavelet based expansions have numerous advantages over Fourier series (e.g., Kurbanmuradov and Sabelfeld, 2008; Walker, 1997) and often lead to stable computations, see (Phoon, Huang, and Quek, 2004).

However, in many cases numerical simulation results need to be confirmed by theoretical analysis. Recently, a considerable attention was given
to the properties of the wavelet transform and of the wavelet orthonormal series representation of
random processes. More information on convergence of wavelet expansions of random processes in various spaces, references and numerous applications can be found in (Cambanis and  Masry, 1994; Didier and  Pipiras, 2008; Istas, 1992; Kozachenko, Olenko and Polosmak, 2010; Kurbanmuradov and Sabelfeld, 2008; Zhang and Waiter, 1994).

In the paper we consider  stationary Gaussian random processes $\mathbf{X}(t)$ and their approximations by sums of wavelet functions
\begin{equation}\label{Xn}\mathbf{X}_{n,\mathbf{k}_n}(t):=\sum_{|k|\le k_0'}\xi_{0k}\phi_{0k}(t)+\sum_{j=0}^{n-1}\sum_{|k|\le k_j}\eta_{jk}\psi_{jk}(t)\,,
\end{equation}
where $\mathbf{k}_n:=(k_0',k_0,...,k_{n-1}).$

Contrary to many theoretical results (see, for example, Kozachenko and Turchyn, 2008; Kurbanmuradov and Sabelfeld, 2008) with infinite series form of $\mathbf{X}_{n,\mathbf{k}_n}(t),$ in direct numerical implementations we always consider truncated series like (\ref{Xn}), where the number of terms in the sums is finite by application reasons.  Most known results (see, for example, Cambanis and  Masry, 1994; Didier and  Pipiras, 2008; Istas, 1992; Wong, 1993; Zhang and Waiter, 1994) concern the mean-square convergence, but for various practical applications one needs to require uniform convergence.

It was shown in (Kozachenko, Olenko and Polosmak, 2010) that, under suitable conditions, the sequence $\mathbf{X}_{n,\mathbf{k}_n}(t)$ converges in probability in Banach space $C([0,T])$, i.e.
\begin{equation}\label{old}P\left\{\sup_{0\le t\le T} |\mathbf{X}(t)-\mathbf{X}_{n,\mathbf{k}_n}(t)|>\varepsilon \right\}\to 0,
\end{equation}
when  $n\to\infty,$ $k_0'\to\infty$ and $k_j\to\infty$ for all $j\in \mathbb{N}_0:=\{0,1,...\}\,.$ More details on random processes in Banach spaces can be found in (Buldygin and  Kozachenko, 2000).

The rate of convergence in (\ref{old}) is another natural question and it is very useful in various computational applications, especially if we are interested in the optimality of the stochastic approximation or the simulations. This question has not been addressed yet.

Our focus in this paper is on fine convergence properties of the wavelet expansions of $\mathbf{X}(t).$
We will obtain an explicit estimate of the rate of convergence in (\ref{old}).
This problem has been posed to the authors several times since the uniform stochastic approximation methods were first proposed. We present the first result on the rate of stochastic uniform convergence of general finite wavelet expansions in the open literature.

We will prove the exponential rapidity of convergence of a wide class of wavelet expansions.
The numbers $n$ and $\mathbf{k}_n$ of terms in the truncated series $\mathbf{X}_{n,\mathbf{k}_n}(t)$ can approach infinity in any arbitrary way. Thought about this way, one sees that the paper deals with the most general class of such wavelet expansions in comparison with particular cases considered by different authors, see, for example, (Cambanis and  Masry, 1994; Kurbanmuradov and Sabelfeld, 2008).

The organization of this article is the following. In the
second section we introduce the necessary background from wavelet
theory and a theorem on uniform convergence in probability of the wavelet expansions of stationary Gaussian random processes, obtained in (Kozachenko, Olenko and Polosmak, 2010).
In \S 3 we specify some upper bounds for the supremum of Gaussian processes and  formulate the main theorem on the rate of convergence. The next section contains the proof of the theorem. Conclusions are made in section 5.

\section{Wavelet representation of random processes}
Let $\phi(x),$ $x\in\mathbb R$ be a function from the space
$L_2(\mathbb R)$ such that $\widehat{\phi}(0)\ne 0$ and  $\widehat{\phi}(y)$ is
continuous at $0,$ where
$$\widehat{\phi}(y)=\int_{\mathbb
R}e^{-iyx}{\phi(x)}\,dx$$
is the Fourier transform of $\phi.$

 Suppose that the following assumption
holds true:
$$\sum_{k\in\mathbb Z} |\widehat{\phi}(y+2{\pi}k)|^2=1\  {\rm (a.e.)}
$$

There exists a function $m_0(x)\in L_2([0,2\pi])$, such that $m_0(x)$
has the period $2\pi$ and
$$\widehat{\phi}(y)=m_0\left(y/2\right)\widehat{\phi}\left(y/2\right)\ {\rm (a.e.)}
$$
 In this case the  function $\phi(x)$ is called
the $f$-wavelet.

Let $\psi(x)$ be the inverse Fourier transform of the function
$$\widehat{\psi}(y)=\overline{m_0\left(\frac
y2+\pi\right)}\cdot\exp\left\{-i\frac
y2\right\}\cdot\widehat{\phi}\left(\frac y2\right).$$
Then the function
$$\psi(x)=\frac1{2\pi}\int_{\mathbb
R}e^{iyx}{\widehat{\psi}(y)}\,dy$$ is called the $m$-wavelet.

Let
\begin{equation}\label{2phijk}\phi_{jk}(x)=2^{j/2}\phi(2^jx-k),\quad
\psi_{jk}(x)=2^{j/2}\psi(2^jx-k),\quad j,k \in\mathbb Z\,.
\end{equation}

It is known that the family of functions $\{\phi_{0k};
\psi_{jk},j\in \mathbb N_0\}$ is an orthonormal basis in
$L_2(\mathbb R)$ (see, for example, (Chui, 1992; Daubechies, 1992; Meyer, 1995)).

An arbitrary function $f(x)\in L_2(\mathbb R)$ can be represented in the form
\begin{equation}\label{2.5}f(x)=\sum_{k\in\mathbb Z}\alpha_{0k}\phi_{0k}(x)+\sum_{j=0}^{\infty}\sum_{k\in\mathbb Z}\beta_{jk}\psi_{jk}(x)\,,
\end{equation}
$$\alpha_{0k}=\int_{\mathbb R}f(x)\overline{\phi_{0k}(x)}\,dx,\quad \beta_{jk}=\int_{\mathbb R}f(x)\overline{\psi_{jk}(x)}\,dx.$$
The representation (\ref{2.5}) is called a wavelet representation.

The series~(\ref{2.5}) converges in the space
$L_2(\mathbb R)$ i.e.
$$\sum_{k\in\mathbb Z}|\alpha_{0k}|^2+\sum_{j=0}^{\infty}\sum_{k\in\mathbb Z}|\beta_{jk}|^2<\infty\,.$$

The integrals $\alpha_{0k}$ and $\beta_{jk}$ may also exist
for functions from $L_p(\mathbb R)$ and other function spaces. Therefore it is possible to obtain the representation~(\ref{2.5}) for function classes which are wider than $L_2(\mathbb R)$ (see, for example, (Jaffard, 2001; Triebel, 2008)).

Let $\{\Omega, \cal{B}, P\}$ be a standard probability space. Let
$\mathbf{X}(t),$ $t\in\mathbb R$ be a  random process
such that $\mathbf E\mathbf{X}(t)=0\,$ for all $t\in\mathbb R$.

It is possible to obtain representations like~(\ref{2.5}) for random
processes, if their sample trajectories are  in the space $L_2(\mathbb R).$ However  the majority of random
processes do not posses this property. For example, sample pathes of
stationary processes are not in  $L_2(\mathbb R)$ (a.s.).

We investigate a representation of the kind~(\ref{2.5}) for
$\mathbf{X}(t)$ with mean-square integrals
$$\xi_{0k}=\int_{\mathbb R}\mathbf{X}(t)\overline{\phi_{0k}(t)}\,dt,\quad \eta_{jk}=\int_{\mathbb R}\mathbf{X}(t)\overline{\psi_{jk}(t)}\,dt\,.$$

Consider the approximants of $\mathbf{X}(t)$ defined by (\ref{Xn}).

\noindent {\bf Assumption S.} (Hardle et al., 1998)  For the function $\phi$ there exists a decreasing function
$\Phi(x),$ $x\ge 0$ such that $\Phi(0)<\infty,$ $|\phi(x)|\le \Phi(|x|)$ (a.e.) and $\displaystyle\int_{\mathbb R}\Phi(|x|)\,dx<\infty\,.$

Let $\mathbf{X}(t)$ be a stationary separable centered Gaussian random process such that
 its covariance function  $R(t,s)=R(t-s)$ is continuous. Let the $f$-wavelet $\phi$
and the corresponding $m$-wavelet $\psi$ be continuous functions and  the assumption {\rm S}  hold true
for both  $\phi$ and  $\psi.$

Theorem~\ref{thuni} below guarantees the uniform convergence of $\mathbf{X}_{n,\mathbf{k}_n}(t)$  to $\mathbf{X}(t).$\vspace{1mm}

\begin{thm}\label{thuni}{\rm(Kozachenko, Olenko and Polosmak, 2010)}  Suppose that the following
conditions hold:
{\rm \begin{enumerate}
  \item\label{con1} there exist $\phi'(u),$ $\widehat{\psi}''(u),$ and  $\widehat{\psi}(0)=0,$ $\widehat{\psi}'(0)=0 ;$
  \item $c_{\phi}:=\sup\limits_{u\in \mathbb R}|\widehat{\phi}(u)|<\infty,$ $
c_{\phi'}:=\sup\limits_{u\in \mathbb R}|\widehat{\phi}'(u)|<\infty,$ $\widehat{\psi}'(u)\in L^1(\mathbb{R}),$ $c_{\psi''}:=\sup\limits_{u\in \mathbb R}|\widehat{\psi}''(u)|<\infty;$
  \item $\widehat{\phi}(u)\to 0$ and $\widehat{\psi}(u)\to 0$ when $u\to \pm\infty;$
  \item\label{con4} there exist $0<\gamma<\frac{1}{2}$  and $\alpha>\frac{1}{2}$ such that $\int\limits_{\mathbb
R}\left(\ln(1+|u|)\right)^{\alpha}|\widehat{\psi}(u)|^{\gamma}\,du<\infty,$\\
 $\int\limits_{\mathbb R}\left(\ln(1+|u|)\right)^{\alpha}|\widehat{\phi}(u)|^{\gamma}\,du<\infty\,;$
  \item\label{con5} there exists $\widehat R(z)$ and  $\sup\limits_{z\in \mathbb R}\widehat R(z)<\infty\,;$
  \item\label{con6} $\int\limits_{\mathbb R}\left|\widehat
R'(z)\right|\,dz<\infty$ and $\int\limits_{\mathbb R}\left|\widehat
R^{(p)}(z)\right||z|^4\,dz<\infty$ for $p=0,1\,.$
\end{enumerate}}
\noindent Then $\mathbf{X}_{n,\mathbf{k}_n}(t)\to \mathbf{X}(t)$ uniformly in probability on each interval $[0,T]$ when  $n\to\infty,$ $k_0'\to\infty$ and $k_j\to\infty$ for all
$j\in \mathbb N_0\,.$
\end{thm}
Before stating the main result, we clarify the role of some of the assumptions.  Two kinds of assumptions are made:
 conditions \ref{con1}-\ref{con4} on the wavelet basis and
 conditions \ref{con5} and \ref{con6} on the random process.

Conditions \ref{con1}-\ref{con4} are related to the smoothness and the decay rate of the wavelet basis functions $\phi$ and $\psi.$ It is easy to check that numerous wavelets satisfy these conditions, for example, the well known Daubechies, Battle-Lemarie and Meyer wavelet bases. Conditions \ref{con5} and \ref{con6} on the random process $\mathbf{X}(t)$ are formulated in terms of the spectral density $\widehat R(z).$ These conditions are related to the behavior of the high-frequency part of the spectrum.

Both sets of assumptions are standard in the convergence studies. On the contrary to many other results in literature, the assumptions are very simple and can be easily verified.

\section{Convergence rate in the space $C[0, T].$}

First we need to specify an estimate for the supremum of Gaussian processes that will be used in the main theorem.

\begin{defin}{\rm(Buldygin and  Kozachenko, 2000, \S3.2)} A set $\mathcal{Q}\subset \mathcal{S} \subset \mathbb{R}$ is called an $\varepsilon$-net in the set $\mathcal{S}$ with respect to the semimetric $\rho$ if for any point $x\in \mathcal{S}$ there exists at least one point $y\in \mathcal{Q}$ such that $\rho(x,y)\le \varepsilon.$
\end{defin}

\begin{defin}{\rm(Buldygin and  Kozachenko, 2000, \S3.2)} Suppose that
$$H_\rho(\mathcal{S},\varepsilon)=\left\{
                          \begin{array}{ll}
                            \ln(N_\rho(\mathcal{S},\varepsilon)), & \hbox{if}\ N_\rho(\mathcal{S},\varepsilon)<+\infty; \\
                            +\infty, & \hbox{if}\ N_\rho(\mathcal{S},\varepsilon)=+\infty,
                          \end{array}
                        \right.
 $$
where $N_\rho(\mathcal{S},\varepsilon)$ is the number of point in a minimal $\varepsilon$-net in the set $\mathcal{S}.$

The function $H_\rho(\mathcal{S},\varepsilon),$ $\varepsilon>0$ is called the metric entropy of the set $\mathcal{S}.$
\end{defin}

\begin{lem}\label{71}{\rm(Buldygin and  Kozachenko, 2000, \S3, (4.10))} Let $Y(t),$ $t\in[0,T]$ be a separable Gaussian random process,
$$\varepsilon_0:=\sup_{t\in[0, T]}\left(E|Y(t)|^2\right)^{1/2}<\infty\,,$$
\begin{equation}\label{76}I(\varepsilon_0):=\frac1{\sqrt2}\int\limits_{0}^{\varepsilon_0}
\left(H(\varepsilon)\right)^{1/2}\,d\varepsilon<\infty\,,\end{equation}
where $H(\varepsilon)$ is the metric entropy of the space $([0, T],\rho),$ $\rho(t,s)=(E|Y(t)-Y(s)|^2)^{1/2}.$

Then
$$P\left\{\sup_{t\in[0,T]} |Y(t)|>u \right\}\le
2\exp\left\{-\frac{\left(u-\sqrt{8uI(\varepsilon_0)}\right)^2}{2\varepsilon_0^2}\right\}\,,$$
where $u>8I(\varepsilon_0).$
\end{lem}

Assume that there exists a nonnegative monotone nondecreasing function
$\sigma(\varepsilon),$ $\varepsilon>0,$ such that $\sigma(\varepsilon)\to 0$ when $\varepsilon\to 0$
and
\begin{equation}\label{74}\sup\limits_{\scriptsize\begin{array}{c}
                     |t-s|\le \varepsilon\\
                      t,s\in [0, T]
                   \end{array}}\left(E|Y(t)-Y(s)|^2\right)^{1/2}\le\sigma(\varepsilon)\,.
\end{equation}
\begin{lem}\label{lem2} Let
\begin{equation}\label{sigma}\sigma(\varepsilon)=\frac c{\left(\ln\left(e^{\alpha}+
\frac1{\varepsilon}\right)\right)^{\alpha}},  {\ } \alpha>1/2\,.\end{equation}
Then {\rm(\ref{76})} holds true and
$$I(\varepsilon_0)\le \delta(\varepsilon_0):=\frac{\gamma}{\sqrt2}\left(\sqrt{\ln(T+1)}+\left({1-\frac{1}{2\alpha}}\right)^{-1}\left(\frac{c}{\gamma}\right)^{\frac{1}{2\alpha}}\right),$$
where
$\gamma:=\min\left(\varepsilon_0,\,\sigma\left(\frac{T}{2}\right)\right)\,.$
\end{lem}
\noindent{\em Proof.}
By simple properties of the metric entropy, see (Buldygin and  Kozachenko, 2000),
$$H(\varepsilon)=\ln N(\varepsilon)\le \ln\left(\frac{T}{2\sigma^{(-1)}(\varepsilon)}+1\right)\,,$$ where
$N(\varepsilon)$ is the number of points in a minimal $\varepsilon$-net in the interval $[0,T].$
Note, that $N(\varepsilon)=1$ if $\varepsilon>\sigma\left(\frac{T}{2}\right).$  Hence $H(\varepsilon)=0$ for $\varepsilon>\sigma\left(\frac{T}{2}\right).$

Hence the integral in (\ref{76}) has the following upper bound

$$I(\varepsilon_0)\le \frac1{\sqrt2}\int\limits_{0}^{\gamma}
\left(\ln\left(\frac{T}{2\sigma^{(-1)}(\varepsilon)}+1\right)\right)^{1/2}\,d\varepsilon\,.$$

The inverse function for $\sigma(\varepsilon)$ is
\begin{equation}\label{inv}
\sigma^{(-1)}(t)=\frac 1{e^{\left(\frac
ct\right)^{1/\alpha}}-e^{\alpha}}\,.\end{equation}
$\sigma^{(-1)}(t)$ is an increasing function for  $t\in\left(0,\,\sigma\left(\frac{T}{2}\right)\right].$ This implies that
$$\frac{T}{2\sigma^{(-1)}(t)}\ge \frac{T}{2\sigma^{(-1)}\left(\sigma\left(\frac{T}{2}\right)\right)}=1.$$

Hence
$$I(\varepsilon_0)\le \frac1{\sqrt2}\int\limits_{0}^{\gamma}
\left(\ln\left(\frac{T}{\sigma^{(-1)}(\varepsilon)}\right)\right)^{1/2}\,d\varepsilon\,.$$

By (\ref{inv}) we get
$$I(\varepsilon_0)= \frac1{\sqrt2}\int\limits_{0}^{\gamma}
\left(\ln\left(T\left(e^{\left(\frac
c\varepsilon\right)^{1/\alpha}}-e^{\alpha}\right)\right)\right)^{1/2}\,d\varepsilon\le \frac1{\sqrt2}\int\limits_{0}^{\gamma}
\left(\ln(T+1)+\left(\frac
c\varepsilon\right)^{1/\alpha}\right)^{1/2}\,d\varepsilon$$
$$\le \frac1{\sqrt2}\int\limits_{0}^{\gamma}
\left(\sqrt{\ln(T+1)}+\left(\frac
c\varepsilon\right)^{\frac{1}{2\alpha}}\right)\,d\varepsilon\le \frac1{\sqrt2}\left(\gamma
\sqrt{\ln(T+1)}+\frac{
c^{\frac{1}{2\alpha}}\gamma^{1-{\frac{1}{2\alpha}}}}{1-{\frac{1}{2\alpha}}}\right)\,.\vspace{-8mm}$$
\hfill $\Box$

\

Now we are ready to formulate the main result of this paper.

\begin{thm} Let $X(t), t\in [0,T]$ be a separable Gaussian stationary random process. Let assumptions {\rm 1-6} of Theorem~{\rm\ref{thuni}} hold true for $X(t).$

Then
$$P\left\{\sup_{t\in[0, T]} |\mathbf{X}(t)-\mathbf{X}_{n,\mathbf{k}_n}(t)|>u \right\}\le
2\exp\left\{-\frac{(u-\sqrt{8u\delta(\varepsilon_{\mathbf{k}_n})})^2}{2\varepsilon_{\mathbf{k}_n}^2}\right\}\,,$$
where $u>8\delta(\varepsilon_{\mathbf{k}_n}),$
$$\varepsilon_{\mathbf{k}_n}:= \sum_{j=0}^{n-1}\frac{A}{2^{j/2}\sqrt{k_j}}+\frac{B}{\sqrt{k_0'}}+\frac{C}{2^{n/2}}.$$
$A,$ $B,$ and $C$ are constants which depend only on the covariance function of $\mathbf{X}(t)$ and the wavelet basis. Explicit expressions for $A,$ $B,$ and $C$ are given in the proof of the theorem.
\end{thm}

\begin{rmk}
$\varepsilon_{\mathbf{k}_n}\to 0,$ if and only if $n,$ $k_0',$ and all $k_j, j\ge 0,$ approach infinity.
\end{rmk}

\begin{rmk}
If $\varepsilon_{\mathbf{k}_n}\to 0$ then $\delta(\varepsilon_{\mathbf{k}_n})\to 0.$ Therefore the convergence in the theorem is exponential with the rate bounded by $\displaystyle 2\exp\left\{-\frac{{const}}{\varepsilon_{\mathbf{k}_n}^2}\right\}.$
\end{rmk}

\begin{rmk}
If we narrow our general class of wavelet expansions $\mathbf{X}_{n,\mathbf{k}_n}(t)$ and impose some additional constraints on rates of the sequences $\mathbf{k}_n$ we can enlarge classes of wavelets bases and random processes in the theorem.
\end{rmk}

\section{Proof of the main theorem}

Let us check that $Y(t)=\mathbf{X}(t)-\mathbf{X}_{n,\mathbf{k}_n}(t)$ satisfies (\ref{74}) with the function $\sigma(\varepsilon)$ given by (\ref{sigma}). One can easily see that
$$\left(E|Y(t)-Y(s)|^2\right)^{1/2}=\left(E\left|\sum_{|k|>k_0'}\xi_{0k}(\phi_{0k}(t)-\phi_{0k}(s))\right.\right.$$
$$\left.\left.+\sum_{j=0}^{n-1}\sum_{|k|>k_j}\eta_{jk}(\psi_{jk}(t)-\psi_{jk}(s))+\sum_{j=n}^{\infty}\sum_{k\in\mathbb
Z}\eta_{jk}(\psi_{jk}(t)-\psi_{jk}(s))\right|
^2\right)^{1/2}$$
$$\le\left(E\left|\sum_{|k|>k_0'}\xi_{0k}(\phi_{0k}(t)-\phi_{0k}(s))\right|
^2\right)^{1/2}+
\sum_{j=0}^{n-1}\left(E\left|\sum_{|k|>k_j}\eta_{jk}(\psi_{jk}(t)-\psi_{jk}(s))\right| ^2\right)^{1/2}$$
$$+\sum_{j=n}^{\infty}\left(E\left|\sum_{k\in\mathbb Z}\eta_{jk}(\psi_{jk}(t)-\psi_{jk}(s))\right| ^2\right)^{1/2}:=\sqrt{S}+\sum_{j=0}^{n-1}\sqrt{S_j}+\sum_{j=n}^\infty \sqrt{R_j}\,.$$

We will give a suitable upper bound for $S_j$, then similar techniques can be used
to deal with the remaining terms $S$ and $R_j.$

$S_j$ satisfies the inequality
 $$S_j\le \sum\limits_{|k|> k_j}\sum\limits_{|l|> k_j}|\mathbf E\eta_{jk}\overline{\eta_{jl}}||\psi_{jk}(t)-\psi_{jk}(s)||\psi_{jl}(t)-\psi_{jl}(s)|\,.$$
Let us consider $\mathbf E\eta_{jk}\overline{\eta_{jl}}.$ By means of Parseval's theorem
we deduce
$${\  }\mathbf E\eta_{jk}\overline{\eta_{jl}}=\int\limits_{\mathbb
R}\int\limits_{\mathbb
R}\mathbf E \mathbf{X}(u)\overline{\mathbf{X}(v)}\ \overline{\psi_{jk}(u)}\psi_{jl}(v)\,dudv=\int\limits_{\mathbb R}\int\limits_{\mathbb R}R(u-v)\overline{\psi_{jk}(u)}\,du\,\psi_{jl}(v)\,dv$$
$$=\int\limits_{\mathbb R}\frac{1}{2 \pi}\int\limits_{\mathbb R}\widehat{R}(z)e^{-ivz}\overline{\widehat{\psi}_{jk}(z)}\,dz\psi_{jl}(v)\,dv=\frac{1}{2 \pi}\int\limits_{\mathbb R}\widehat{R}(z)\,\overline{\widehat{\psi}_{jk}(z)}\,\widehat{\psi}_{jl}(z)\,dz\,.$$

The order of integration can be changed because
$$ \int\limits_{\mathbb R}\int\limits_{\mathbb R}\left|\widehat{R}(z)e^{-ivz}\overline{\widehat{\psi}_{jk}(z)}\psi_{jl}(v)\right|\,dz\,dv\le \sup_{z\in\mathbb R}|\widehat{R}(z)|\cdot
\int\limits_{\mathbb R}|{\widehat{\psi}_{jk}(z)}|\,dz\cdot\int\limits_{\mathbb R}|\psi_{jl}(v)|\,dv<\infty\,.$$
The last expression is finite due to  (\ref{2phijk}), the assumption S, the estimate (\ref{c2}), and the representation
\begin{equation}\label{104.11}\widehat{\psi}_{jk}(z)=\frac{e^{-i\frac k{2^j}z}}{2^{j/2}}\cdot\widehat{\psi}\left(\frac z{2^j}\right).
\end{equation}

We begin with the case $k\not= l,$ $k\not= 0,$ $l\not= 0.$

Applying integration by parts,
the theorem's assumptions, and the inequality $\left|\widehat{\psi}(\frac
z{2^j})\right|\le c_{\psi''}\frac {|z|^2}{2^{2j}}$ yields the following
 $$|\mathbf E\eta_{jk}\overline{\eta_{jl}}|=\left|\frac{1}{2 \pi}\int\limits_{\mathbb R}\widehat{R}(z)
 \frac{e^{i\frac {k-l}{2^j}z}}{2^{j}}\left|\widehat{\psi}\left(\frac z{2^j}\right)\right|^2\,dz\right|=\left|\frac{1}{\pi 2^{j+1}}\right.\cdot \frac{2^j}{i(k-l)}
\left[\left.\widehat{R}(z)
 e^{i\frac {k-l}{2^j}z}\left|\widehat{\psi}\left(\frac z{2^j}\right)\right|^2\right|_{z=-\infty}^{+\infty}\right.$$
 $$-\int\limits_{\mathbb R}\left(\widehat{R}'(z)\left|\widehat{\psi}\left(\frac z{2^j}\right)\right|^2
 \left.\left.+\widehat{R}(z)\frac{2}{2^j}\,\Re\left(\overline{\widehat{\psi}\left(\frac z{2^j}\right)}
\widehat{\psi}'\left(\frac z{2^j}\right)\right)\right)e^{i\frac
{k-l}{2^j}z}\,dz\right]\right|$$
\begin{equation}\label{aphi}\le\frac{1}{2 \pi|k-l|}\int\limits_{\mathbb R}\left(\left|\widehat{R}'(z)\right|
\left|\widehat{\psi}\left(\frac
z{2^j}\right)\right|^2+\frac{2}{2^j}\left|\widehat{R}(z)\right|\,\left|\overline{\widehat{\psi}\left(\frac
z{2^j}\right)}\widehat{\psi}'\left(\frac z{2^j}\right)\right|\right)\,dz\le\frac{A^{\psi}}{2^{4j}|k-l|}\,,
\end{equation}
where $$A^{\psi}:=\frac{c_{\psi''}^2}{2 \pi}\int\limits_{\mathbb
R}\left(|\widehat{R}'(z)||z|^4
+2|\widehat{R}(z)||z|^3\right)\,dz<\infty\,,$$
because theorem's conditions \ref{con5} and \ref{con6}.

We use (\ref{104.11}) to estimate the term $|{\psi}_{jl}(t)-{\psi}_{jl}(s)|\,.$ Then
$${\psi}_{jl}(t)=\int\limits_{\mathbb R}\frac {e^{itz}e^{-i\frac l{2^j}z}}{\pi2^{j/2+1}}\widehat{\psi}\left(\frac z{2^j}\right)\,dz=\int\limits_{\mathbb R}\frac{e^{it\left(z+\frac{2^j}l\pi\right)}e^{-i\left(\frac l{2^j}z+\pi\right)}}{\pi 2^{j/2+1}}\widehat{\psi}\left(\frac z{2^j}+\frac {\pi}l\right)\,dz\,.
$$
Hence
$${\psi}_{jl}(t)=\frac{1}{2^{j/2+2}\pi}\int\limits_{\mathbb R}e^{-i\frac l{2^j}z}\left(e^{itz}\widehat{\psi}\left(\frac z{2^j}\right)-e^{it\left(z+\frac{2^j}l\pi\right)}\widehat{\psi}\left(\frac z{2^j}+\frac {\pi}l\right)\right)\,dz\,.
$$
Therefore
$$|{\psi}_{jl}(t)-{\psi}_{jl}(s)|=\frac{1}{2^{j/2+2}\pi}\left|\int\limits_{\mathbb R}e^{-i\frac
l{2^j}z}\left[\left(e^{itz}\widehat{\psi}\left(\frac
z{2^j}\right)-e^{it\left(z+\frac{2^j}l\pi\right)}\widehat{\psi}\left(\frac
z{2^j}+\frac {\pi}l\right)\right)\right.\right.$$
$$-\left.\left.\left(e^{isz}\widehat{\psi}\left(\frac
z{2^j}\right)-e^{is\left(z+\frac{2^j}l\pi\right)}\widehat{\psi}\left(\frac
z{2^j}+\frac {\pi}l\right)\right)\right]\,dz\right|=\left|\frac z{2^j}=u\right|\le\frac{2^{j/2-2}}{\pi}\int\limits_{\mathbb
R}\left(\left|e^{it2^ju}-e^{is2^ju}\right.\right.$$
\begin{equation}\label{I2}\left.\left.-e^{it2^j\left(u+\frac{\pi}l\right)}+e^{is2^j\left(u+\frac{\pi}l\right)}
\right|\left|\widehat{\psi}(u)\right|+\left|e^{it2^j\left(u+\frac{\pi}l\right)}-e^{is2^j\left(u+\frac{\pi}l\right)}
\right|\left|\widehat{\psi}(u)-\widehat{\psi}\left(u+\frac
{\pi}l\right)\right|\right)\,du\,.\end{equation}
By the inequality (59) given in (Kozachenko and Rozora, 2003) for $\alpha>0$ we obtain
\begin{equation}\label{105}|e^{itz}-e^{isz}|=\left(2-2\cos\left(z(t-s)\right)\right)^{1/2}=2\left|\sin\left(\frac{z(t-s)}2\right)\right|
\le2\left(\frac{\ln\left(e^{\alpha}+ \frac{|z|}2\right)}{\ln\left(e^{\alpha}+
\frac1{|t-s|}\right)}\right)^{\alpha} \,.
\end{equation}
An application of this inequality to the second part of the integral in (\ref{I2}) and the substitution $u+\frac {\pi}l=v$ result~in
$$I_2=\int\limits_{\mathbb
R}\left|e^{it2^jv}-e^{is2^jv} \right|\left|\widehat{\psi}\left(v-\frac
{\pi}l\right)-\widehat{\psi}(v)\right|\,dv\le\frac2{\left(\ln\left(e^{\alpha}+
\frac1{|t-s|}\right)\right)^{\alpha}}$$
$$\times\int\limits_{\mathbb
R}\left(\ln\left(e^{\alpha}+
\frac{2^j|v|}2\right)\right)^{\alpha}\left|\widehat{\psi}\left(v-\frac
{\pi}l\right)-\widehat{\psi}(v)\right|^{\beta}\left|\widehat{\psi}\left(v-\frac
{\pi}l\right)-\widehat{\psi}(v)\right|^{1-\beta}\,dv,$$
where $\beta:=1-\gamma\in (1/2,1).$

In the following derivations, we will use assumptions \ref{con1}, \ref{con4}
and the estimates
$$\left|\widehat{\psi}\left(v-\frac{\pi}l\right)-\widehat{\psi}(v)\right|^{\beta}\le c_{\psi'}^\beta\left(\frac{\pi}l\right)^\beta,$$
$$\left|\widehat{\psi}\left(v-\frac
{\pi}l\right)-\widehat{\psi}(v)\right|^{1-\beta}\le 2^{1-\beta}\left(\left|\widehat{\psi}\left(v-\frac
{\pi}l\right)\right|^{1-\beta}+\left|\widehat{\psi}(v)\right|^{1-\beta}\right)\,,$$
$$\int\limits_{\mathbb
R}\left(\ln\left(e^{\alpha}+
\frac{2^j|v|}2\right)\right)^{\alpha}\left|\widehat{\psi}\left(v\right)\right|^{1-\beta}\,dv\le\int\limits_{\mathbb
R}\left(\ln\left[5^{j+1}\left(\frac{e^{\alpha}}{5^{j+1}}+ \frac{2^{j-1}}{5^{j+1}}\left|v\right|\right)\right]\right)^{\alpha}$$
\begin{equation}\label{5}\times\left| \widehat{\psi}(v)\right|^{1-\beta}\,dv\le 2^\alpha\left((\ln5)^\alpha (j+1)^{\alpha}c_0+c_1\right)<\infty\,,\end{equation}
$$\int\limits_{\mathbb
R}\left(\ln\left(e^{\alpha}+
\frac{2^j|v|}2\right)\right)^{\alpha}\left|\widehat{\psi}\left(v-\frac
{\pi}l\right)\right|^{1-\beta}\,dv\le\int\limits_{\mathbb
R}\left(\ln\left[5^{j+1}\left(\frac{e^{\alpha}}{5^{j+1}}+ \frac{2^{j-1}}{5^{j+1}}\left(\left|v\right|+\frac{\pi}{|l|}\right)\right)\right]\right)^{\alpha}$$
$$\times \left|\widehat{\psi}(v)\right|^{1-\beta}\,dv
\le 2^\alpha\left((\ln5)^\alpha (j+1)^{\alpha}c_0+c_1\right)<\infty\,,$$
where
$$c_0:=\int\limits_{\mathbb
R}\left|\widehat{\psi}(v)\right|^{1-\beta}\,dv<\infty\,,\quad c_1:=\int\limits_{\mathbb
R}\left(\ln(1+|v|)\right)^{\alpha}\left|\widehat{\psi}(v)\right|^{1-\beta}\,dv<\infty\,.$$

The integral $c_0$ is finite because of the boundedness of $\widehat{\psi}(v)$ and condition~\ref{con4}.

Using these facts, we get
  \begin{equation}\label{205}I_2\le\frac{2^{3+\alpha-\beta}\pi^\beta c_{\psi'}^\beta}{|l|^{\beta}\left(\ln\left(e^{\alpha}+
\frac1{|t-s|}\right)\right)^{\alpha}}\left((\ln5)^\alpha (j+1)^{\alpha}c_0+c_1\right)\,.\end{equation}

 We can similarly estimate the first part of the integral. It is easy to see that
$$|\Delta|:\le\left|e^{it2^ju}-e^{is2^ju}\right|\left|1-e^{it2^j\frac{\pi}l}\right|+
\left|e^{is2^j\frac {\pi}l}-e^{it2^j\frac{\pi}l}\right|$$
$$\le2\left|\sin\left(\frac{2^ju(t-s)}2\right)\right|\cdot\left|\sin\left(\frac{2^j\pi t}{2l}\right)\right|+2\left|\sin\left(\frac{2^j\pi(t-s)}{2l}\right)\right|\,.$$
Note that, by (\ref{105}) and (\ref{5}):
\begin{equation}\label{+}\left|\sin\left(\frac{2^ju(t-s)}2\right)\right|\le\left(\frac{\ln(e^{\alpha}+ \frac{2^j|u|}2)}{\ln\left(e^{\alpha}+
\frac1{|t-s|}\right)}\right)^{\alpha}
\le 2^\alpha\frac{(j+1)^\alpha (\ln 5)^\alpha+\left(\ln\left(1+ |u|\right)\right)^\alpha}{\left(\ln\left(e^{\alpha}+
\frac1{|t-s|}\right)\right)^{\alpha}}\,.\end{equation}
Due to (Kozachenko,  E.  Turchyn, 2008, Lemma 4.2)
\begin{equation}\label{*}\left|\sin\left(\frac{2^j\pi(t-s)}{2l}\right)\right|\le \left|\frac{2^j\pi(t-s)}{2l}\right|\le \frac{2^jc_{\alpha}}{|l|\left(\ln\left(e^{\alpha}+ \frac1{|t-s|}\right)\right)^{\alpha}},\quad  \alpha>0\,,\end{equation}
where $c_{\alpha}$ depends only on $T $ and $\alpha\,.$

Applying inequalities (\ref{+}), (\ref{*}) and
$$\left|\sin\left(\frac{2^j\pi t}{2l}\right)\right|\le \left|\frac{2^j\pi t}{2l}\right|
\le \frac{2^{j-1}\pi T}{|l|}$$
 we get
$$|\Delta|\le\frac{2^{j+1}\left(\pi T 2^{\alpha-1}\left[(j+1)^\alpha (\ln 5)^\alpha+\left(\ln\left(1+ |u|\right)\right)^\alpha\right]+c_{\alpha}\right)}{|l|\left(\ln\left(e^{\alpha}+ \frac1{|t-s|}\right)\right)^{\alpha}}
\,.$$

Using these inequalities, the first part of the integral can be estimated as following
$$I_1=\int\limits_{\mathbb
R}\left|e^{it2^ju}-e^{is2^ju}-e^{it2^j\left(u+\frac{\pi}l\right)}+e^{is2^j\left(u+\frac{\pi}l\right)}
\right|\left|\widehat{\psi}(u)\right|du $$
\begin{equation}\label{206}\le\frac{2^{j+1}}{|l|\left(\ln\left(e^{\alpha}+
\frac1{|t-s|}\right)\right)^{\alpha}}\cdot \left(\pi T 2^{\alpha-1}\left((j+1)^\alpha (\ln 5)^\alpha c_2+c_3\right)+c_{\alpha}c_2\right)\,,\end{equation}
where
$$c_2:=\int\limits_{\mathbb
R}\left|\widehat{\psi}(v)\right|\,dv<\infty\,,\quad c_3:=\int\limits_{\mathbb
R}\left(\ln(1+|v|)\right)^{\alpha}\left|\widehat{\psi}(v)\right|\,dv<\infty\,.$$

The integrals $c_2$ and $c_3$ are finite because $\widehat{\psi}(v)$ is bounded:
\begin{equation}\label{c2}c_2\le\sup\limits_{u\in\mathbb R}\left|\widehat{\psi}(u)\right|\int\limits_{\mathbb
R}\frac{\left|\widehat{\psi}(v)\right|^{1-\beta}}{\left(\sup_{u\in\mathbb R}\left|\widehat{\psi}(u)\right|\right)^{1-\beta}}\,dv=
\left(\sup\limits_{u\in\mathbb R}\left|\widehat{\psi}(u)\right|\right)^{\beta}c_0<\infty\,,
\end{equation}
$$c_3\le \left(\sup\limits_{u\in\mathbb R}\left|\widehat{\psi}(u)\right|\right)^{\beta}c_1<\infty\,.$$

Using (\ref{205}) and (\ref{206}), we obtain:
$$|{\psi}_{jk}(t)-{\psi}_{jk}(s)|\cdot|{\psi}_{jl}(t)-{\psi}_{jl}(s)|\le
\frac{2^{j-4}}{\pi^2 |k|^{\beta}|l|^{\beta}\left(\ln\left(e^{\alpha}+
\frac1{|t-s|}\right)\right)^{2\alpha}}$$
$$\times\left(2^{3+\alpha-\beta}\pi^\beta c_{\psi'}^\beta\left((\ln5)^\alpha (j+1)^{\alpha}c_0+c_1\right)+
2^{j+1}\left(\pi T 2^{\alpha-1}\left((j+1)^\alpha  \right.\right.\right.$$
\begin{equation}\label{ff}\left.\left.\left. \times(\ln 5)^\alpha c_2+c_3\right)+c_{\alpha}c_2\right) \right)^2\le
\frac{(j+1)^{2\alpha}2^{3j-2}K^2}{|k|^{\beta}|l|^{\beta}\left(\ln\left(e^{\alpha}+
\frac1{|t-s|}\right)\right)^{2\alpha}}\,,
\end{equation}
where
$$K:= \pi^{-1}\left(2^{3+\alpha-\beta}\pi^\beta c_{\psi'}^\beta\left((\ln5)^\alpha c_0+c_1\right)+\pi T 2^{\alpha-1}\left((\ln 5)^\alpha c_2+c_3\right)+c_{\alpha}c_2\right)\,.$$

Thus
\begin{equation}\label{knel}\sum\limits_{\scriptsize\begin{array}{c}
                     |k|> k_j,|l|> k_j\\
                      k\ne l
                   \end{array}}
                   \hspace{-5mm}|\mathbf E\eta_{jk}\overline{\eta_{jl}}||\psi_{jk}(t)-\psi_{jk}(s)|
|\psi_{jl}(t)-\psi_{jl}(s)|\le\frac{(j+1)^{2\alpha}A^{\psi}QK^2}{2^{j}\left(\ln\left(e^{\alpha}+
\frac1{|t-s|}\right)\right)^{2\alpha}} \,,
\end{equation}
where
$$Q:=
\sum\limits_{k> k_j} \sum\limits_{l> k_j}\frac{1}{2(l+k)k^{\beta}l^{\beta}}+\sum\limits_{k>l} \sum\limits_{
l>k_j}\frac {1}{(k-l)l^{\beta}k^{\beta}}\le\sum\limits_{k> k_j} \sum\limits_{l> k_j}\frac {1}{4\sqrt {kl}k^{\beta}l^{\beta}} +\sum\limits_{m=1}^\infty \sum\limits_{l>
k_j}\frac {1}{ml^{\beta}(m+l)^{\beta}}\,.$$

Then, using the inequality $x+y>c_{\delta}^{-1}x^{\delta}y^{1-\delta}, {\ }x,y>0, {\
}0<\delta<1,$ where $c_{\delta}:={\delta^{\delta}(1-\delta)^{1-\delta}}$, we obtain
\begin{equation}\label{star}Q\le\left(\sum\limits_{k> k_j} \frac {1}{2k^{\frac12+\beta}}\right)^2
+c^\beta_{\delta}\sum\limits_{m=1}^\infty\frac {1}{m^{1+\delta\beta}} \sum\limits_{l> k_j}\frac
{1}{l^{(2-\delta)\beta}}\end{equation}
$$\le\left(\sum\limits_{k=1}^\infty \frac {1}{2k^{\frac12+\beta}}\right)^2
+c^\beta_{\delta}\sum\limits_{m=1}^\infty\frac {1}{m^{1+\delta\beta}} \sum\limits_{l=1}^\infty\frac
{1}{l^{(2-\delta)\beta}}=:Q_1<\infty\,.$$
 The statement becomes apparent if $\delta\in (0,2-1/\beta)$ is chosen, as $\beta \in \left(\frac{1}{2},1\right).$

 We can similarly exploit the case $l\ne0, k=0:$
$$|\mathbf E\eta_{j0}\overline{\eta_{jl}}|=\left|\frac{1}{2^{j+1}\pi}\right.\left(-\frac{2^j}{il}\right)\left[\left.\widehat R(z)
e^{-i\frac l{2^j}z}\left|\widehat{\psi}\left(\frac
z{2^j}\right)\right|^2\right|_{z=-\infty}^{+\infty}\right.$$
$$-\int\limits_{\mathbb R}\left(\widehat{R}'(z)\left|\widehat{\psi}\left(\frac z{2^j}\right)\right|^2
\left.\left.+\widehat R(z)\frac{2}{2^j}\,\Re\left(\widehat{\psi}\left(\frac
z{2^j}\right)\widehat{\psi}'\left(\frac z{2^j}\right)\right)\right) e^{-i\frac
l{2^j}z}\,dz\right]\right|$$
$$\le\frac{1}{2 \pi|l|}\int\limits_{\mathbb R}\left(\left|\widehat{R}'(z)\right|
\left|\widehat{\psi}\left(\frac
z{2^j}\right)\right|^2+\frac{2}{2^j}\left|\widehat{R}(z)\right|\,\left|\widehat{\psi}\left(\frac
z{2^j}\right) \widehat{\psi}'\left(\frac z{2^j}\right)\right|\right)\,dz$$
$$\le\frac{1}{2 \pi|l|}\int\limits_{\mathbb R}\left(|\widehat{R}'(z)|\left(c_{\psi''}\frac {|z|^2}{2^{2j}}\right)^2
+2|\widehat{R}(z)|c_{\psi''}\frac {|z|^2}{2^{2j}}\frac {|z|}{2^{j}}
c_{\psi''}\right)\,dz=\frac{A^{\psi}}{2^{4j}|l|}\,.
$$

By (\ref{105}) and (\ref{5})
$$|{\psi}_{j0}(t)-{\psi}_{j0}(s)|=\frac{2^{j/2}}{2
\pi}\left|\int\limits_{\mathbb R}(e^{it2^{j}u}-e^{is2^{j}u})
\widehat{\psi}(u)\,du\right|$$
$$\le\frac{2^{j/2}}{
\pi}\int\limits_{\mathbb
R}\left(\frac{\ln\left(e^{\alpha}+\frac{2^j|u|}2\right)}{\ln\left(e^{\alpha}+
\frac1{|t-s|}\right)}\right)^{\alpha}|\widehat{\psi}(u)|\,du\le\frac{2^{j/2}}{
\pi\left(\ln\left(e^{\alpha}+
\frac1{|t-s|}\right)\right)^{\alpha}}$$
\begin{equation}\label{0}\times\int\limits_{\mathbb
R}\left((j+1)\ln 5+\ln(1+|u|)\right)^{\alpha}|\widehat{\psi}(u)|\,du
\le\frac{2^{j/2+\alpha}\left((j+1)^\alpha\left(\ln 5\right)^\alpha c_2+c_3\right)}{\pi\left(\ln\left(e^{\alpha}+ \frac1{|t-s|}\right)\right)^{\alpha}}\,.\end{equation}

Then
$$\sum\limits_{|l|> k_j}|\mathbf E\eta_{j0}\overline{\eta_{jl}}||\psi_{j0}(t)-\psi_{j0}(s)||\psi_{jl}(t)-\psi_{jl}(s)|\le 2\sum\limits_{l=1}^\infty\frac{A^{\psi}}{2^{4j}l}\cdot
\frac{2^{j/2+\alpha}(j+1)^\alpha\left(\left(\ln 5\right)^\alpha c_2+c_3\right)}{\pi\left(\ln\left(e^{\alpha}+ \frac1{|t-s|}\right)\right)^{\alpha}}$$
\begin{equation}\label{105.1}\times  \frac
{2^{3j/2-1}(j+1)^{\alpha}K}{l^{\beta}\left(\ln\left(e^{\alpha}+
\frac1{|t-s|}\right)\right)^{\alpha}}
\le\frac{(j+1)^{2\alpha}\cdot q}{2^{2j}\left(\ln\left(e^{\alpha}+
\frac1{|t-s|}\right)\right)^{2\alpha}}\,,
\end{equation}
where
$$q:=\frac{2^{\alpha}A^{\psi}K((\ln 5)^\alpha c_2+c_3)}{\pi}\cdot\sum\limits_{l=1}^\infty
\frac{1}{l^{1+\beta}}<\infty\,.$$
Similarly
\begin{equation}\label{105.2}\sum\limits_{|k|> k_j}|\mathbf E\eta_{jk}\overline{\eta_{j0}}||\psi_{jk}(t)-\psi_{jk}(s)||\psi_{j0}(t)-\psi_{j0}(s)|
\le\frac{(j+1)^{2\alpha}\cdot q}{2^{2j}\left(\ln\left(e^{\alpha}+
\frac1{|t-s|}\right)\right)^{2\alpha}}\,.
\end{equation}

For the case $l=k, k\not=0$ we get the estimates
$$|{\psi}_{jk}(t)-{\psi}_{jk}(s)|^2\le
\frac{(j+1)^{2\alpha}2^{3j-2}K^2}{|k|^{2\beta}\left(\ln\left(e^{\alpha}+
\frac1{|t-s|}\right)\right)^{2\alpha}}\,,$$
\begin{equation}\label{kl}\mathbf E|\eta_{jk}|^2=\frac{1}{2 \pi}\int\limits_{\mathbb R}\widehat{R}(z)\,
 2^{-j}\left|\widehat{\psi}\left(\frac z{2^j}\right)\right|^2\,dz\le\frac{1}{2^{j+1} \pi}
\int\limits_{\mathbb R}|\widehat{R}(z)|\left(c_{\psi''}\frac {|z|^2}{2^{2j}}\right)^2\,dz=\frac{A^{\psi}_1}{2^{5j}}\,,\end{equation}
where $$A^{\psi}_1:=\frac{c_{\psi''}^2}{2 \pi}\int\limits_{\mathbb
R}|\widehat{R}(z)||z|^4\,dz<\infty\,.$$

Hence
$$\sum\limits_
{\scriptsize\begin{array}{c}
                     |k|> k_j\\
                       l=k
                   \end{array}}
\mathbf |E\eta_{jk}\overline{\eta_{jl}}|\cdot |\psi_{jk}(t)-\psi_{jk}(s)||\psi_{jl}(t)-\psi_{jl}(s)|= 2\sum\limits_{k=k_j+1}^{\infty} \mathbf E|\eta_{jk}|^2|{\psi}_{jk}(t)-{\psi}_{jk}(s)|^2$$
\begin{equation}\label{kel}\le
\sum\limits_{k=k_j+1}^{\infty} \frac{(j+1)^{2\alpha} A^{\psi}_1 K^2}{2^{2j+1}k^{2\beta}\left(\ln\left(e^{\alpha}+
\frac1{|t-s|}\right)\right)^{2\alpha}}\le \frac{(j+1)^{2\alpha}\cdot q_1}{2^{2j}\left(\ln\left(e^{\alpha}+
\frac1{|t-s|}\right)\right)^{2\alpha}}\,,
\end{equation}
where $$q_1:= \frac{A^{\psi}_1 K^2}{2}\cdot\sum\limits_{k=1}^{\infty}\frac{1}{k^{2\beta}}<\infty.$$

Finally, for $k=l=0,$ applying (\ref{kl}) and (\ref{0}),  we get
\begin{equation}\label{kl0}
\mathbf E|\eta_{j0}|^2\cdot|{\psi}_{j0}(t)-{\psi}_{j0}(s)|^2\le \frac{(j+1)^{2\alpha}\cdot q_2}{2^{4j}\left(\ln\left(e^{\alpha}+ \frac1{|t-s|}\right)\right)^{2\alpha}}\,, \end{equation}
where
$$q_2:=\frac{2^{2\alpha}A^{\psi}_1}{\pi^2} \left(\left(\ln 5\right)^\alpha c_2+c_3\right)^2\,.$$

It follows from (\ref{knel}), (\ref{star}) and (\ref{kel}) that
$$\sum\limits_{j=0}^{n-1}\sqrt{S_j}=\sum\limits_{j=0}^{n-1}\left(\mathbf E\left|\sum\limits_{|k|> k_j}\eta_{jk}(\psi_{jk}(t)-\psi_{jk}(s))\right|^2\right)^{1/2}\le\sum\limits_{j=0}^{n-1}\left(\frac{(j+1)^{2\alpha}\cdot q_1}{2^{2j}\left(\ln\left(e^{\alpha}+
\frac1{|t-s|}\right)\right)^{2\alpha}}
\right.$$
$$\left.+
\frac{(j+1)^{2\alpha}A^{\psi}Q_1K^2}
{2^{j}\left(\ln\left(e^{\alpha}+ \frac1{|t-s|}\right)\right)^{2\alpha}}
\right)^{1/2}
\le \frac{B_0}{\left(\ln\left(e^{\alpha}+
\frac1{|t-s|}\right)\right)^{\alpha}}\,,
$$
where
$$B_0:=\left(
q_1+A^{\psi}Q_1K^2\right)^{1/2}\cdot\sum\limits_{j=0}^{\infty}\frac{(j+1)^{\alpha}}{2^{j/2}}<\infty\,.$$

Using (\ref{knel}), (\ref{star}), (\ref{105.1}), (\ref{105.2}), (\ref{kel}) and (\ref{kl0}) we obtain
$$\sum_{j=n}^\infty \sqrt{R_j}\le\sum_{j=n}^\infty \left( \frac{(j+1)^{2\alpha}A^{\psi}Q_1K^2}{2^{j}\left(\ln\left(e^{\alpha}+
\frac1{|t-s|}\right)\right)^{2\alpha}}+\frac{2(j+1)^{2\alpha}\cdot q}{2^{2j}\left(\ln\left(e^{\alpha}+
\frac1{|t-s|}\right)\right)^{2\alpha}}\right.$$
$$\left.+\frac{(j+1)^{2\alpha}\cdot q_1}{2^{2j}\left(\ln\left(e^{\alpha}+
\frac1{|t-s|}\right)\right)^{2\alpha}}+\frac{(j+1)^{2\alpha}\cdot q_2}{2^{4j}\left(\ln\left(e^{\alpha}+ \frac1{|t-s|}\right)\right)^{2\alpha}}\right)^{1/2}\le \frac{B_1}{\left(\ln\left(e^{\alpha}+ \frac1{|t-s|}\right)\right)^{\alpha}}\,, $$
where
$$B_1:=\left(
q+q_1+q_2+A^{\psi}QK^2\right)^{1/2}\cdot\sum\limits_{j=0}^{\infty}\frac{(j+1)^{\alpha}}{2^{j/2}}<\infty\,.$$

The analysis for $S$ is similar. First we evaluate  $|\mathbf E\xi_{0k}\overline{\xi_{0l}}|$ in the case $k\not= l\,:$
$$|\mathbf E\xi_{0k}\overline{\xi_{0l}}|=\left|\frac{1}{2 \pi}\int\limits_{\mathbb R}\widehat R(z)\overline{\widehat{\phi}_{0k}(z)}\,
\widehat{\phi}_{0l}(z)\,dz\right|=\left|\frac{1}{2
\pi}\int\limits_{\mathbb R}\widehat
R(z)e^{i(k-l)z}\left|\widehat{\phi}( z)\right|^2 \,dz\right|$$
$$\le\frac{1}{2 \pi|k-l|}\,\left|\left.\widehat R(z)e^{i(k-l)z}\right.|\widehat{\phi}(z)|^2
\right|_{z=-\infty}^{+\infty}-\int\limits_{\mathbb
R}\left.\left(\widehat R'(z)|\widehat{\phi}(z)|^2 +\widehat
R(z)2\Re\left(\overline{\widehat{\phi}(z)}\widehat{\phi}'(z)\right)\right)\right.$$
$$\times\left.e^{i(k-l)z}\,dz\right| \le\frac{1}{2 \pi|k-l|}\cdot\int\limits_{\mathbb
R}\left(\left|\widehat R'(z)\right||\widehat{\phi}(z)|^2 +2|\widehat
R(z)|\, \left|\widehat{\phi}(z)\widehat{\phi}'(z)\right|\right)\,dz\le\frac{ A^{\phi}}{|k-l|}\,,
$$
where
$$A^{\phi}:=\frac{1}{2 \pi}\,\left(c_{\phi}^2\int\limits_{\mathbb
R}\left|\widehat R'(z)\right|\,dz
+2c_{\phi}c_{\phi'}\int\limits_{\mathbb R}\left|\widehat
R(z)\right|\,dz\right)<\infty\,$$
because of  assumptions \ref{con5} and \ref{con6}.

Similarly to (\ref{ff}) we  derive
$$|\phi_{0k}(t)-\phi_{0k}(s)|\cdot|\phi_{0l}(t)-\phi_{0l}(s)|
\le \frac {(K^{\phi})^2}{4|l|^{\beta}|k|^{\beta}\left(\ln\left(e^{\alpha}+
\frac1{|t-s|}\right)\right)^{2\alpha}}\,,
$$ where
$$K^{\phi}:=\pi^{-1}\left(2^{3+\alpha-\beta}\pi^\beta c_{\phi'}^\beta\left((\ln5)^\alpha c_{\phi 0}+c_{\phi 1}\right)+\pi T 2^{\alpha-1}\left((\ln 5)^\alpha c_{\phi 2}+c_{\phi 3}\right)+c_{\alpha}c_{\phi 2}\right)\,,$$
$$c_{\phi 0}:=\int\limits_{\mathbb
R}\left|\widehat{\phi}(v)\right|^{1-\beta}\,dv<\infty\,,\quad c_{\phi 1}:=\int\limits_{\mathbb
R}\left(\ln(1+|v|)\right)^{\alpha}\left|\widehat{\phi }(v)\right|^{1-\beta}\,dv<\infty\,,$$
$$c_{\phi  2}:=\int\limits_{\mathbb
R}\left|\widehat{\phi }(v)\right|\,dv<\infty\,,\quad c_{\phi  3}:=\int\limits_{\mathbb
R}\left(\ln(1+|v|)\right)^{\alpha}\left|\widehat{\phi }(v)\right|\,dv<\infty\,.$$

Therefore
$$\sum\limits_{\scriptsize\begin{array}{c}
                     |k|> k_0',|l|> k_0'\\
                      k\ne l
                   \end{array}}|
\mathbf E\xi_{0k}\overline{\xi_{0l}}||\phi_{0k}(t)-\phi_{0k}(s)|
|\phi_{0l}(t)-\phi_{0l}(s)|$$
\begin{equation}\label{1070}\le \sum\limits_{\scriptsize\begin{array}{c}
                     |k|> k_0',|l|> k_0'\\
                      k\ne l
                   \end{array}}
\frac{A^{\phi}}{|k-l|}\cdot \frac
{(K^{\phi})^2}{4|l|^{\beta}|k|^{\beta}\left(\ln\left(e^{\alpha}+
\frac1{|t-s|}\right)\right)^{2\alpha}}
\le
\frac {A^{\phi}(K^{\phi})^2 Q}{\left(\ln\left(e^{\alpha}+
\frac1{|t-s|}\right)\right)^{2\alpha}}\,.
\end{equation}

For the case $l=k,$ $k\not= 0$ we obtain
$$|\phi_{0k}(t)-\phi_{0k}(s)|^2\le \frac {(K^{\phi})^2}{4\,|k|^{2\beta}\left(\ln\left(e^{\alpha}+
\frac1{|t-s|}\right)\right)^{2\alpha}}\,,$$
$$\mathbf E|\xi_{0k}|^2=\frac{1}{2 \pi}\int\limits_{\mathbb
R}\widehat R(z)|\widehat{\phi}( z)|^2dz\,\le \frac{c_\phi^2}{2\pi}\int\limits_{\mathbb
R}|\widehat R(z)|\,dz=:A_1^\phi<\infty\,.
$$

Hence
$$\sum\limits_
{\scriptsize\begin{array}{c}
                     |k|> k_0'\\
                       l=k
                   \end{array}}
\mathbf |E\xi_{0k}\overline{\xi_{0l}}|\cdot |\phi_{0k}(t)-\phi_{0k}(s)||\phi_{0l}(t)-\phi_{0l}(s)|= 2\sum\limits_{k=k_0'+1}^{\infty} \mathbf E|\xi_{0k}|^2|{\phi}_{0k}(t)-{\phi}_{0k}(s)|^2$$
\begin{equation}\label{kel1}\le
\sum\limits_{k=k_0'+1}^{\infty} \frac{ A^{\phi}_1 (K^\phi)^2}{2k^{2\beta}\left(\ln\left(e^{\alpha}+
\frac1{|t-s|}\right)\right)^{2\alpha}}< \frac{q_{\phi 1}}{\left(\ln\left(e^{\alpha}+
\frac1{|t-s|}\right)\right)^{2\alpha}}\,,
\end{equation}
where $$q_{\phi 1}:= \frac{ A^{\phi}_1 (K^\phi)^2}{2}\cdot\sum\limits_{k=1}^{\infty}\frac{1}{k^{2\beta}}<\infty.$$

Combining (\ref{1070}) with (\ref{kel1}) one gets that
$$\sqrt{S}=\left(\mathbf E\left|\sum\limits_{|k|> k_0'}\xi_{0k}(\phi_{0k}(t)-\phi_{0k}(s))\right|^2\right)^{1/2}
\le \frac {B_2}{\left(\ln\left(e^{\alpha}+
\frac1{|t-s|}\right)\right)^{\alpha}},
$$
where
$$B_2:=\left(q_{\phi 1}+A^{\phi}(K^{\phi})^2 Q\right)^{1/2}\,.$$

Finally combining all estimates, we obtain
$$\left(E\left|(\mathbf{X}(t)-\mathbf{X}_{n,\mathbf{k}_n}(t))-(\mathbf{X}(s)-\mathbf{X}_{n,\mathbf{k}_n}(s))\right|^2\right)^{1/2}
\le \frac{B_0+B_1+B_2}{\left(\ln(e^{\alpha}+
\frac1{h})\right)^{\alpha}}=:\sigma(h)\,.
$$

Then by lemmata~\ref{71} and \ref{lem2}
$$P\left\{\sup_{t\in[0,T]} |\mathbf{X}(t)-\mathbf{X}_{n,\mathbf{k}_n}(t)|>u \right\}\le
2\exp\left\{-\frac{(u-\sqrt{8uI(\tilde{\varepsilon}_{\mathbf{k}_n})})^2}{2\tilde{\varepsilon}_{\mathbf{k}_n}^2}\right\}\,,$$
where $$\tilde{\varepsilon}_{\mathbf{k}_n}=\sup_{t\in[0,T]} \left(E\left|\mathbf{X}(t)-\mathbf{X}_{n,\mathbf{k}_n}(t)\right|^2\right)^{1/2}\,.$$

Let us investigate $\tilde{\varepsilon}_{\mathbf{k}_n}$ as a function of $\mathbf{k}_n.$ One can easily see that
$$\left(E\left|\mathbf{X}(t)-\mathbf{X}_{n,\mathbf{k}_n}(t)\right|^2\right)^{1/2}\le\left(E\left|\sum_{|k|>k_0'}\xi_{0k}\phi_{0k}(t)\right| ^2\right)^{1/2}+$$
$$+\sum_{j=0}^{n-1}\left(E\left|\sum_{|k|>k_j}\eta_{jk}\psi_{jk}(t)\right| ^2\right)^{1/2}+\sum_{j=n}^{\infty}\left(E\left|\sum_{k\in\mathbb Z}\eta_{jk}\psi_{jk}(t)\right| ^2\right)^{1/2}\,.$$
Let us consider the second sum
$$\sum_{j=0}^{n-1}\left(E\left|\sum_{|k|>k_j}\eta_{jk}\psi_{jk}(t)\right| ^2\right)^{1/2}\le\sum_{j=0}^{n-1}\left(\sum_{|k|>k_j}\sum_{|l|>k_j}E|\eta_{jk}\eta_{jl}||\psi_{jk}(t)||\psi_{jl}(t)|\right)^{1/2}\,.$$
For $l\not=0$ one can estimate $\psi_{jl}$ as follows:
$$|{\psi}_{jl}(t)|=\frac {2^{j/2}}{2 \pi|l|}\left|\left(e^{itz}\widehat{\psi}\left(\frac z{2^j}\right)\left.e^{-i\frac l{2^j}z}\right)\right|_{z=-\infty}^{+\infty}-\int\limits_{\mathbb R}\left(\frac{e^{itz}}{2^j}\widehat{\psi}'\left(\frac z{2^j}\right)+ite^{itz}\widehat{\psi}\left(\frac z{2^j}\right)\right)e^{-i\frac l{2^j}z}\,dz\right|$$
\begin{equation}\label{81}\le\frac {2^{j/2}}{2 \pi|l|}\left(\int\limits_{\mathbb R}\left|\widehat{\psi}'(u)\right|\,du+2^jT\int\limits_{\mathbb R}\left|\widehat{\psi}(u)\right|\,du\right)\le\frac {2^{3j/2}}{|l|}B^{\psi}_1\,,
\end{equation}
where
$$B^{\psi}_1:=\frac{1}{2 \pi}\left(\int\limits_{\mathbb R}\left|\widehat{\psi}'(u)\right|\,du+T\int\limits_{\mathbb R}\left|\widehat{\psi}(u)\right|\,du\right)<\infty\,.$$
Then by (\ref{aphi}), (\ref{81}) and (\ref{star}) with $\beta=1, \delta=1/2$ we obtain:
$$\sum_{j=0}^{n-1}\left(E\left|\sum_{|k|>k_j}\eta_{jk}\psi_{jk}(t)\right| ^2\right)^{\frac12}\le
\sum_{j=0}^{n-1}\left(\sum\limits_{\scriptsize\begin{array}{c}
                     |k|> k_j,|l|> k_j\\
                      k\ne l
                   \end{array}}\frac {A^{\psi}\left(B^{\psi}_1\right)^2}{2^{j}|k-l||k||l|}+\sum\limits_{|k|> k_j}\frac {A_1^{\psi}\left(B^{\psi}_1\right)^2}{2^{2j}|k|^2}
\right)^{\frac12}$$
$$\le B^{\psi}_1\sum_{j=0}^{n-1}\left(\frac{4A^{\psi}}{2^j} \left(\left(\sum\limits_{k> k_j} \frac {1}{2k^{\frac32}}\right)^2
+\frac{1}{2}\sum\limits_{m=1}^\infty\frac {1}{m^{3/2}} \sum\limits_{l> k_j}\frac
{1}{l^{3/2}} \right)+\frac{2A_1^{\psi}}{2^{2j}}\sum_{k>k_j}^\infty\frac{1}{k^2}\right)^{1/2}\le  \sum_{j=0}^{n-1}\frac{A}{2^{j/2}k_j^{1/2}},$$
where
$$A:=B^{\psi}_1\left(6 A^{\psi}\sum\limits_{m=1}^\infty\frac {1}{m^{3/2}}+4 A_1^{\psi}\right)^{1/2}\,.$$

In the above introduced transformations we used the inequality
$$\sum_{k>k_j}\frac 1{k^{\alpha}}< \sum_{k=k_j}^{\infty}\int\limits_{k}^{k+1}\frac{dx}{x^{\alpha}}=\int\limits_{k_j}^{\infty}\frac{dx}{x^{\alpha}}=\frac{1}{(\alpha-1)k_j^{\alpha-1}},\quad \alpha>1.$$

Similarly to (\ref{81}),  for $j=0$ we get $|\phi_{0k}(t)|\le {B^{\phi}_1}/{|k|},
$ where
$$B^{\phi}_1:=\frac{1}{2
\pi}\left(\int\limits_{\mathbb
R}\left|\widehat{\phi}'(u)\right|\,du+T\int\limits_{\mathbb
R}\left|\widehat{\phi}(u)\right|\,du\right)<\infty\,.$$

Therefore,  by (\ref{star}) with $\beta=1, \delta=1/2$ we obtain
$$\left(E\left|\sum_{|k|>k_0'}\xi_{0k}\phi_{0k}(t)\right| ^2\right)^{1/2}\le\Biggl(\sum\limits_{\scriptsize\begin{array}{c}
                     |k|> k_0',|l|> k_0'\\
                      k\ne l
                   \end{array}}\frac {A^{\phi}\left(B^{\phi}_1\right)^2}{|k-l||k||l|}+\sum\limits_{|k|> k_0'}\frac {A_1^{\psi}\left(B^{\psi}_1\right)^2}{|k|^2}
\Biggr)^{1/2}\le\frac{B}{\sqrt{k_0'}} ,$$
where
$$B:=B^{\phi}_1\left(6A^{\phi}\sum\limits_{m=1}^\infty\frac {1}{m^{3/2}}+4A_1^{\phi} \right)^{1/2}.$$

$\psi_{j0}(t)$ can be bounded as follows:
\begin{equation}\label{psi0}|{\psi}_{j0}(t)|=\left|\frac{1}{2 \pi}\int\limits_{\mathbb R}\frac {e^{itz}}{2^{j/2}}\widehat{\psi}\left(\frac z{2^j}\right)\,dz\right|\le
\frac{2^{j/2-1}}{\pi}\int\limits_{\mathbb R}|\widehat{\psi}(z)|dz=\frac{2^{j/2-1}}{\pi}\,c_2\,.\end{equation}

By (\ref{aphi}), (\ref{81}), (\ref{psi0}) and (\ref{star}) with $\beta=1, \delta=1/2$ we obtain:
$$\sum_{j=n}^{\infty}\left(E\left|\sum_{k\in\mathbb Z}\eta_{jk}\psi_{jk}(t)\right| ^2\right)^{1/2}\le \sum_{j=n}^{\infty}\Biggl(\sum\limits_{\scriptsize\begin{array}{c}
                     k,l\in\mathbb Z\\
                      k\ne l, kl\ne 0
                   \end{array}}\frac {A^{\psi}\left(B^{\psi}_1\right)^2}{2^j|k-l||k||l|}+\sum\limits_{\scriptsize\begin{array}{c}
                     k\in\mathbb Z\\
                      k\ne 0
                   \end{array}}\frac {A_1^{\psi}\left(B^{\psi}_1\right)^2}{2^{2j} |k|^2}$$
                   $$+\sum\limits_{\scriptsize\begin{array}{c}
                     k\in\mathbb Z\\
                      k\ne 0
                   \end{array}}\frac {c_2 A^{\psi}B_1^{\psi}}{\pi 2^{2j} |k|^2}+\frac {c_2^2 A_1^{\psi}}{4\pi^2 2^{4j}}\Biggr)^{1/2}\le\sum_{j=n}^{\infty}\frac{1}{2^{j/2}}\cdot
\left(4A^{\psi} \left(B^{\psi}_1\right)^2 \left(\left(\sum\limits_{k=1} \frac {1}{2k^{\frac32}}\right)^2
+\frac{1}{2}\left(\sum\limits_{m=1}^\infty\frac {1}{m^{3/2}}\right)^2 \right)\right.$$
$$\left. +\sum_{k=1}^\infty\frac{1}{k^2}\left(A_1^{\psi} \left(B^{\psi}_1\right)^2+\frac{c_2 A^{\psi}B_1^{\psi}}{\pi}\right)+\frac {c_2^2 A_1^{\psi}}{32\pi^2}\right)^{1/2}\le \frac{C}{2^{n/2}}, $$
where
$$C:=(2+\sqrt{2})\left(3A^{\psi} \left(B^{\psi}_1\right)^2 \left(\sum\limits_{k=1} \frac {1}{k^{\frac32}}\right)^2 +\left(A_1^{\psi} \left(B^{\psi}_1\right)^2+\frac{c_2 A^{\psi}B_1^{\psi}}{\pi}\right)\sum_{k=1}^\infty\frac{1}{k^2}+\frac {c_2^2 A_1^{\psi}}{32\pi^2}\right)^{1/2}.$$

Therefore the upper bound for $\tilde{\varepsilon}_{\mathbf{k}_n}$ is determined as follows $$\tilde{\varepsilon}_{\mathbf{k}_n}=\sup_{t\in[0,T]}\left(E|\mathbf{X}(t)-\mathbf{X}_{n,\mathbf{k}_n}(t)|^2\right)^{1/2}
\le \sum_{j=0}^{n-1}\frac{A}{2^{j/2}\sqrt{k_j}}+\frac{B}{\sqrt{k_0'}}+\frac{C}{2^{n/2}}=\varepsilon_{\mathbf{k}_n}.$$

It is worth noticing that $\delta(\cdot)$ is an increasing function, $I(\varepsilon_0)\le \delta(\varepsilon_0)$ for any $\varepsilon_0,$ and $\tilde{\varepsilon}_{\mathbf{k}_n}\le \varepsilon_{\mathbf{k}_n}.$ Hence for $u>8\delta(\varepsilon_{\mathbf{k}_n})$ we get

$$\exp\left\{-\frac{\left(u-\sqrt{8uI(\tilde{\varepsilon}_{\mathbf{k}_n})}\right)^2}{2\tilde{\varepsilon}_{\mathbf{k}_n}^2}\right\}\le \exp\left\{-\frac{\left(u-\sqrt{8u\delta(\varepsilon_{\mathbf{k}_n})}\right)^2}{2\varepsilon_{\mathbf{k}_n}^2}\right\}\,.$$

Finally, an application of Lemmata~\ref{71} and \ref{lem2} completes the proof.\hfill $\Box$

\section{Conclusions}
We have obtained the rate of uniform convergence of wavelet
expansions of stationary Gaussian random processes.  The wavelet expansions are more general than the Fourier-wavelets decompositions studied in (Kurbanmuradov and Sabelfeld, 2008). The most general form of the expansions was studied.
The results are obtained under simple conditions which can be easily verified. The conditions  are weaker than those in the former literature. The obtained theorem is the first result on the rate of uniform convergence of general finite wavelet expansions for stochastic processes in the open literature.

This analysis is new and provides a constructive algorithm for determining the number of terms in the wavelet expansions to ensure the uniform approximation of stochastic processes with given accuracy.

It should be mentioned that in the general case the random coefficients in (\ref{Xn}) may be dependent and form an overcomplete system of basis functions. That is why the implementation of the proposed method has promising potential for nonstationary random processes, (see, for example, (Meyer, Sellan, and Taqqu, 1999)). However, the generalizations to various nonstationary random processes are not always straightforward, because conditions \ref{con5} and \ref{con6} must be changed. The essence of the problem is a modification of the proof by obtaining new upper bounds on $|\mathbf E\eta_{jk}\overline{\eta_{jl}}|$ and $|\mathbf E\xi_{0k}\overline{\xi_{0l}}|.$ For some nonstationary classes of stochastic processes this non-trivial problem requires a full-length paper itself.

The main aim of this paper was to derive results which are valid for general wavelet expansions. It would be of interest to adopt and specify the results for
different wavelet bases, new classes of stochastic processes, and to examine the tightness of the estimates.

\section{Acknowledgements}
This work was partly supported by
La Trobe University Research Grant-501821 "Sampling, wavelets and optimal stochastic modelling" and La Trobe University Research Grant "Stochastic Approximation in Finance and Signal Processing".

\vskip 3mm

\noindent \textbf{\Large References}
\vskip 3mm

\noindent Buldygin, V.V.,  Kozachenko, Yu.V. (2000). \textit{Metric Characterization of Random Variables and Random Processes.}  Providence R.I.: American Mathematical Society.

\vskip 3mm
\noindent Cambanis, S.,   Masry, E. (1994). Wavelet approximation of deterministic and random signals: convergence properties and rates. \textit{IEEE Trans. Inf. Theory } 40(4):1013--1029.

\vskip 3mm
\noindent Chui, C.K. (1992).  \textit{An Introduction to Wavelets.} New York: Academic Press.

\vskip 3mm
\noindent Daubechies,  I. (1992). \textit{Ten Lectures on Wavelets.} Philadelphia: SIAM.

\vskip 3mm
\noindent Didier, G.,  Pipiras, V. (2008). Gaussian stationary processes: adaptive wavelet decompositions, discrete approximations and their convergence. \textit{J. Fourier Anal. and Appl.} \mbox{14:203--234.}

\vskip 3mm
\noindent Hardle,  W., Kerkyacharian, G., Picard,  D.,  Tsybakov, A. (1998).  \textit{Wavelets, Approximation and Statistical Applications.}  New York: Springer.

\vskip 3mm
\noindent Istas, J.  (1992).  Wavelet coefficients of a gaussian process and applications.  \textit{Annales de l'institut Henri Poincar\'{e} (B) Probabilit\'{e}s et Statistiques}  28(4):537--556.

\vskip 3mm
\noindent Jaffard, S. (2001). Wavelet expansions, function spaces and multifractal analysis. In:  Byrnes, J.S., ed., \textit{Twentieth century harmonic analysis -- a celebration.} Dordrecht: Kluwer Acad. Publ., 127--144.

\vskip 3mm
\noindent Kozachenko, Yu.,  Olenko, A., Polosmak, O. (2011). Uniform convergence of wavelet expansions of Gaussian random processes.
\textit{Stochastic Analysis and Applications} 29(2):169--184.

\vskip 3mm
\noindent Kozachenko, Yu., Rozora, I. (2003).
 Simulation of Gaussian stochastic processes. \textit{Random Operators Stoch. Equations} 11(3):275--296.

\vskip 3mm
\noindent Kozachenko, Yu.,  Turchyn,  E. (2008).  Conditions of uniform convergence of wavelet expansion of
$\phi$-sub-Gaussian random processes. \textit{Theory Prob. Math. Statist.} 78:74--85.

\vskip 3mm
\noindent Kurbanmuradov, O., Sabelfeld, K.  (2008). Convergence of fourier-wavelet models for Gaussian random processes. \textit{SIAM J. Numer. Anal.} 46(6):3084--3112.

\noindent Meyer, Y. (1995). \textit{Wavelets and Operators.} Cambridge: Cambridge University Press.

\vskip 3mm
\noindent Meyer, Y., Sellan, F., Taqqu, M.S. (1999). Wavelets, generalized white noise and fractional integration:
the synthesis of fractional Brownian motion. \textit{J. Fourier Anal. and Appl.} 5(5):465--494.

\vskip 3mm
\noindent Phoon, K.K.,  Huang, H.W.,  Quek, S.T. (2004). Comparison betveen Karhunen-Loeve and
wavelet expansions for simulation of Gaussian processes. \textit{Comput.  Structures} \mbox{82:985–-991.}

\vskip 3mm
\noindent Triebel, H. (2008). \textit{Function Spaces and Wavelets on Domains.} Zu\"{u}rich: European Mathematical Society.

\vskip 3mm
\noindent Wong, P.W. (1993). Wavelet decomposition of harmonizable random processes. \textit{IEEE Trans. Inf. Theory} 39(1):7--18.

\vskip 3mm
\noindent Walker, J. S. (1997). Fourier analysis and wavelet analysis. \textit{Notices Amer. Math. Soc.} \mbox{44: 658--670.}

\vskip 3mm
\noindent Zhang, J., Waiter, G. (1994). A wavelet-based KL-like expansion for wide-sense stationary random processes. \textit{IEEE Trans. Signal Proc.} 42(7):1737--1745.

\end{document}